\documentclass[11pt]{article}
\newtheorem{theorem}{Theorem}[section]
\newtheorem{lemma}[theorem]{Lemma}

\newtheorem{cor}[theorem]{Corollary}

\newtheorem{conjecture}[theorem]{Conjecture}

\usepackage{color}
\usepackage{amssymb,amsmath}
\usepackage{wasysym}
\usepackage{graphicx}

\usepackage{longtable}
\usepackage{array}

\numberwithin{equation}{section}
\def\f{\noindent}
\def\h{\hfill  $\Box$\vspace{10pt}}
\def\p{\f {\bf Proof.}\hskip10pt}

\begin{document}
\title{On the conjecture about the exponential reduced Sombor index}
\author{Wei Gao\\
Department of Mathematics\\ Pennsylvania State University at Abington,\\
Abington, PA, 19001, USA\\
E-mail: wvg5121@psu.edu}
\date{}
\maketitle

\begin{abstract}
Let $G=(V(G),E(G))$ be a graph and $d(v)$ be the degree of the vertex $v\in V(G)$.
The exponential reduced Sombor index of $G$, denoted by $e^{SO_{red}}(G)$, is defined as
$$e^{SO_{red}}(G)=\sum_{uv\in E(G)}e^{\sqrt{(d(u)-1)^2+(d(v)-1)^2}}.$$

We obtain a characterization of extremal trees with the maximal
exponential reduced Sombor index among all chemical trees of order $n$.
This result shows the conjecture on the exponential reduced Sombor index proposed
by Liu, You, Tang and Liu [On the reduced Sombor index and its applications,
MATCH Commun. Math. Comput. Chem. 86 (2021) 729--753] is negative.

\bigskip
\f {\it AMS classification:}  05C50, 05C09, 05C92.

\f {\bf Keywords:} Chemical tree; exponential reduced Sombor index; extremal tree.
\end{abstract}

\section{Introduction}
\hskip\parindent
Let $G$ be a graph with vertex set $V(G)$
and edge set $E(G)$, and $d(v)$ be the degree of the vertex $v\in V(G)$.
A tree $T$ is a chemical tree (or molecular tree) if $d_T(v)\le 4$ for $v\in V(T)$.
Let ${\cal CT}_n$ be the set of all chemical trees of order $n$.

The reduced Sombor index of $G$, denoted by $SO_{red}(G)$, is defined as
$$SO_{red}(G)=\sum_{uv\in E(G)}\sqrt{(d(u)-1)^2+(d(v)-1)^2}.$$
The exponential reduced Sombor index of $G$, denoted by $e^{SO_{red}}(G)$, is defined as
$$e^{SO_{red}}(G)=\sum_{uv\in E(G)}e^{\sqrt{(d(u)-1)^2+(d(v)-1)^2}}.$$
An interesting conjecture on exponential reduced Sombor index was proposed in \cite{match86n3_729-753} as follows.

\begin{conjecture}[\cite{match86n3_729-753}]\label{conj-1}\
Let $T\in {\cal CT}_n$, $n\ge 5$. Then
$$e^{SO_{red}}(T)\le \left\{\begin{array}{ll}
\frac{2}{3}(n+1)e^3+\frac{1}{3}(n-5)e^{3\sqrt{2}}, & n\equiv 2 (\!\!\!\!\mod 3);\\
\frac{1}{3}(2n+1)e^3+\frac{1}{3}(n-13)e^{3\sqrt{2}}+3 e^{\sqrt{13}}, & n\equiv 1 (\!\!\!\!\mod 3);\\
\frac{1}{3}ne^3+\frac{1}{3}(n-9)e^{3\sqrt{2}}+2 e^{\sqrt{10}}, & n\equiv 0 (\!\!\!\!\mod 3).
\end{array}\right.$$
\end{conjecture}

Liu, You, Tang, Liu \cite{match86n3_729-753} proposed two
conjectures on the reduced Sombor index and the exponential reduced Sombor index, respectively.
The paper \cite{match88n3_583-591} solved one conjecture on the reduced Sombor index and
pointed out that the other conjecture is still open.
For more results for the reduced Sombor index and the exponential reduced Sombor index of graphs, we
refer to \cite{match86n3_729-753,match88n3_583-591}.

In this note, we obtain a characterization of extremal trees with the maximal
exponential reduced Sombor index among all chemical trees of order $n$.
This result shows that Conjecture \ref{conj-1} is negative.

\section{Lemmas}
\hskip\parindent
For a chemical tree $T\in {\cal CT}_n$, let $n_i(T)$ be the number of vertices with degree $i$,
and $m_{ij}(T)$ be the number of edges joining a vertex of degree $i$ and a vertex of degree $j$.
We call that the edge joining a vertex of degree $i$ and a vertex of degree $j$ is an $(i,j)$-edge.
It is clear that
\begin{align}
    &n_1(T)+n_2(T)+n_3(T)+n_4(T) = n, \label{eq2.1}\\
    &n_1(T)+2n_2(T)+3n_3(T)+4n_4(T) = 2(n-1), \label{eq2.2}\\
    &m_{12}(T)+m_{13}(T)+m_{14}(T) = n_1(T),  \label{eq2.3}\\
    &m_{12}(T)+2m_{22}(T)+m_{23}(T)+m_{24}(T) = 2n_2(T), \label{eq2.4}\\
    &m_{13}(T)+m_{23}(T)+2m_{33}(T)+m_{34}(T) = 3n_3(T), \label{eq2.5}\\
    &m_{14}(T)+m_{24}(T)+m_{34}(T)+2m_{44}(T) = 4n_4(T). \label{eq2.6}
\end{align}
Moreover, Gutman and Miljkovi\'{c} in \cite{Match41(2000)57-70} established the following relations:
\begin{align}
m_{14}(T)=&\frac{2n+2}{3}-\frac{4}{3}m_{12}(T)-\frac{10}{9}m_{13}(T)-\frac{2}{3}m_{22}(T)-\frac{4}{9}m_{23}(T)\nonumber\\
&-\frac{1}{3}m_{24}(T)-\frac{2}{9}m_{33}(T)-\frac{1}{9}m_{34}(T), \label{eq2.7}
\\
m_{44}(T)=&\frac{n-5}{3}+\frac{1}{3}m_{12}(T)+\frac{1}{9}m_{13}(T)-\frac{1}{3}m_{22}(T)-\frac{5}{9}m_{23}(T)\nonumber\\
&-\frac{2}{3}m_{24}(T)-\frac{7}{9}m_{33}(T)-\frac{8}{9}m_{34}(T). \label{eq2.8}
\end{align}

Denote $f(x,y)=e^{\sqrt{(x-1)^2+(y-1)^2}}$.
Then for $T\in {\cal CT}_n$, by (\ref{eq2.7}) and (\ref{eq2.8}),
\begin{align}
e^{SO_{red}}(T)=&\sum_{1\le i\le j\le 4}f(i,j)m_{ij}(T)\nonumber\\
=&f(1,2)m_{12}(T)+f(1,3)m_{13}(T)+f(1,4)m_{14}(T)+f(2,2)m_{22}(T)+f(2,3)m_{23}(T)\nonumber\\
&+f(2,4)m_{24}(T)+f(3,3)m_{33}(T)+f(3,4)m_{34}(T)+f(4,4)m_{44}(T)\nonumber\\
=&\frac{2n+2}{3}f(1,4)+\frac{n-5}{3}f(4,4)+M_{12}m_{12}(T)+M_{13}m_{13}(T)+M_{22}m_{22}(T)\nonumber\\
&+M_{23}m_{23}(T)+M_{24}m_{24}(T)+M_{33}m_{33}(T)+M_{34}m_{34}(T),\label{eq2.9}
\end{align}
where
\begin{align*}
&M_{12}=f(1,2)-\frac{4}{3}f(1,4)+\frac{1}{3}f(4,4),\\
&M_{13}=f(1,3)-\frac{10}{9}f(1,4)+\frac{1}{9}f(4,4),\\
&M_{22}=f(2,2)-\frac{2}{3}f(1,4)-\frac{1}{3}f(4,4),\\
&M_{23}=f(2,3)-\frac{4}{9}f(1,4)-\frac{5}{9}f(4,4),\\
&M_{24}=f(2,4)-\frac{1}{3}f(1,4)-\frac{2}{3}f(4,4),\\
&M_{33}=f(3,3)-\frac{2}{9}f(1,4)-\frac{7}{9}f(4,4),\\
&M_{34}=f(3,4)-\frac{1}{9}f(1,4)-\frac{8}{9}f(4,4).
\end{align*}

With the help of Mathematical software, we can verify that Lemma \ref{lem-0} holds.

\begin{lemma}\label{lem-0}\
(1)\ $M_{33}<M_{23}<M_{22}<M_{24}<M_{34}<M_{13}<M_{12}<0$;

(2)\ $M_{12}+M_{22}+M_{24}<2M_{12}+2M_{24}<2M_{13}+M_{34}<M_{12}+M_{24}$.
\end{lemma}

\begin{lemma}\label{lem-1}\
Let $n\ge 5$ and $T\in {\cal CT}_n$.
If $n_2(T)=n_3(T)=0$, then $n\equiv 2 (\!\!\!\!\mod 3)$, and
\begin{align}\label{eq2.10}
e^{SO_{red}}(T)=\frac{2n+2}{3}f(1,4)+\frac{n-5}{3}f(4,4).
\end{align}
\end{lemma}

\p Let $n_2(T)=n_3(T)=0$. By (\ref{eq2.1}) and (\ref{eq2.2}),
we have $n_1(T)+n_4(T)=n$ and $n_1(T)+4n_4(T)=2(n-1)$.
Then $n_1(T)=\frac{2n+2}{3}$ and $n_4(T)=\frac{n-2}{3}$.
It implies that $n\equiv 2 (\!\!\!\!\mod 3)$.

Note that $m_{12}(T)=m_{13}(T)=m_{22}(T)=m_{23}(T)=m_{24}(T)=m_{33}(T)=m_{34}(T)=0$.
So by (\ref{eq2.9}),
$$e^{SO_{red}}(T)=\frac{2n+2}{3}f(1,4)+\frac{n-5}{3}f(4,4).$$
The lemma holds. \h

\begin{lemma}\label{lem-2}\
Let $n\ge 5$ and $T\in {\cal CT}_n$.
If $n_2(T)=0$ and $n_3(T)=1$, then $n\equiv 1 (\!\!\!\!\mod 3)$, and
\begin{align}\label{eq2.11}
e^{SO_{red}}(T)\le 2f(1,3)+f(3,4)+\frac{2n-5}{3}f(1,4)+\frac{n-7}{3}f(4,4),
\end{align}
the equality holds if and only if $m_{13}(T)=2$ and $m_{34}(T)=1$.
\end{lemma}

\p Let $n_2(T)=0$ and $n_3(T)=1$. By (\ref{eq2.1}) and (\ref{eq2.2}),
$n_1(T)+n_4(T)=n-1$ and $n_1(T)+4 n_4(T)=2n-5$.
Then $n_1(T)=\frac{2n+1}{3}$ and $n_4(T)=\frac{n-4}{3}$.
It implies that $n\equiv 1 (\!\!\!\!\mod 3)$.

Since $n_2(T)=0$ and $n_3(T)=1$, we have $m_{12}(T)=m_{22}(T)=m_{23}(T)=m_{24}(T)=m_{33}(T)=0$.
Note that $m_{13}(T)+m_{34}(T)=3$, and $0\le m_{13}(T)\le 2$.
By Lemma \ref{lem-0}, $M_{34}<M_{13}<0$.
Then by (\ref{eq2.9}),
\begin{align*}
e^{SO_{red}}(T)=&\frac{2n+2}{3}f(1,4)+\frac{n-5}{3}f(4,4)+M_{13}m_{13}(T)+M_{34}m_{34}(T)
\\
\le &\frac{2n+2}{3}f(1,4)+\frac{n-5}{3}f(4,4)+2M_{13}+M_{34}
\\
=&\frac{2n+2}{3}f(1,4)+\frac{n-5}{3}f(4,4)+2\left(f(1,3)-\frac{10}{9}f(1,4)+\frac{1}{9}f(4,4)\right)\\
&+\left(f(3,4)-\frac{1}{9}f(1,4)-\frac{8}{9}f(4,4)\right)
\\
=&2f(1,3)+f(3,4)+\frac{2n-5}{3}f(1,4)+\frac{n-7}{3}f(4,4).
\end{align*}
It is easy to see that the equality in (\ref{eq2.11}) holds if and only if $m_{13}(T)=2$ and $m_{34}(T)=1$.
\h

\begin{lemma}\label{lem-3}\
Let $n\ge 5$ and $T\in {\cal CT}_n$.
If $n_2(T)=1$ and $n_3(T)=0$, then $n\equiv 0 (\!\!\!\!\mod 3)$, and
\begin{align}\label{eq2.12}
e^{SO_{red}}(T)\le f(1,2)+f(2,4)+\frac{2n-3}{3}f(1,4)+\frac{n-6}{3}f(4,4),
\end{align}
the equality holds if and only if $m_{12}(T)=m_{24}(T)=1$.
\end{lemma}

\p Let $n_2(T)=1$ and $n_3(T)=0$. By (\ref{eq2.1}) and (\ref{eq2.2}),
we have $n_1(T)+n_4(T)=n-1$ and $n_1(T)+4n_4(T)=2n-4$.
Then $n_1(T)=\frac{2n}{3}$ and $n_4(T)=\frac{n-3}{3}$.
It implies that $n\equiv 0 (\!\!\!\!\mod 3)$.

Since $n_2(T)=1$ and $n_3(T)=0$, we have $m_{13}(T)=m_{22}(T)=m_{23}(T)=m_{33}(T)=m_{34}(T)=0$.
Note that $m_{12}(T)+m_{24}(T)=2$ and $0\le m_{12}(T)\le 1$.
By Lemma \ref{lem-0} and (\ref{eq2.9}), we have $M_{24}<M_{12}<0$, and
\begin{align*}
e^{SO_{red}}(T)=&\frac{2n+2}{3}f(1,4)+\frac{n-5}{3}f(4,4)+M_{12}m_{12}(T)+M_{24}m_{24}(T)
\\
\le &\frac{2n+2}{3}f(1,4)+\frac{n-5}{3}f(4,4)+M_{12}+M_{24}
\\
=&\frac{2n+2}{3}f(1,4)+\frac{n-5}{3}f(4,4)+\left(f(1,2)-\frac{4}{3}f(1,4)+\frac{1}{3}f(4,4)\right)\\
&+\left(f(2,4)-\frac{1}{3}f(1,4)-\frac{2}{3}f(4,4)\right)
\\
=&f(1,2)+f(2,4)+\frac{2n-3}{3}f(1,4)+\frac{n-6}{3}f(4,4).
\end{align*}
It is easy to see that the equality in (\ref{eq2.12}) holds if and only if $m_{12}(T)=m_{24}(T)=1$.
\h

\begin{lemma}\label{lem-5}\
Let $n\ge 5$ and $T\in {\cal CT}_n$.
If $n_2(T)+n_3(T)\ge 2$, then
\begin{align}\label{eq2.13}
&e^{SO_{red}}(T)<2f(1,3)+f(3,4)+\frac{2n-5}{3}f(1,4)+\frac{n-7}{3}f(4,4).
\end{align}
\end{lemma}

\p We consider the following three cases.

{\bf Case 1.}\ $n_2(T)=0$ and $n_3(T)\ge 2$.

Note that $m_{12}(T)=m_{22}(T)=m_{23}(T)=m_{24}(T)=0$.
Then by (\ref{eq2.9}),
\begin{align*}
e^{SO_{red}}(T)=&\frac{2n+2}{3}f(1,4)+\frac{n-5}{3}f(4,4)
+M_{13}m_{13}(T)+M_{33}m_{33}(T)+M_{34}m_{34}(T).
\end{align*}
Since $n_3(T)\ge 2$, we have $m_{13}(T)+m_{33}(T)+m_{34}(T)\ge 5$, and $m_{33}(T)+m_{34}(T)\ge 1$.
By Lemma \ref{lem-0}, $M_{33}<M_{34}<M_{13}<0$.
So
\begin{align*}
e^{SO_{red}}(T)\le &\frac{2n+2}{3}f(1,4)+\frac{n-5}{3}f(4,4)+4M_{13}+M_{34}
\\
< &\frac{2n+2}{3}f(1,4)+\frac{n-5}{3}f(4,4)+2M_{13}+M_{34}
\\
=&\frac{2n+2}{3}f(1,4)+\frac{n-5}{3}f(4,4)+2\left(f(1,3)-\frac{10}{9}f(1,4)+\frac{1}{9}f(4,4)\right)\\
&+\left(f(3,4)-\frac{1}{9}f(1,4)-\frac{8}{9}f(4,4)\right)
\\
=&2f(1,3)+f(3,4)+\frac{2n-5}{3}f(1,4)+\frac{n-7}{3}f(4,4).
\end{align*}

{\bf Case 2.}\ $n_2(T)\ge 2$ and $n_3(T)=0$.

Note that $m_{13}(T)=m_{23}(T)=m_{33}(T)=m_{34}(T)=0$.
Then by (\ref{eq2.9}),
\begin{align*}
e^{SO_{red}}(T)=&\frac{2n+2}{3}f(1,4)+\frac{n-5}{3}f(4,4)
+M_{12}m_{12}(T)+M_{22}m_{22}(T)+M_{24}m_{24}(T).
\end{align*}
Let $P=v_1v_2\dots v_t$ be a longest path in $T$ with $d_T(v_1)=d_T(v_t)=2$.
Denote $N_T(v_1)=\{v_0,v_2\}$ and $N_T(v_t)=\{v_{t-1},v_{t+1}\}$.
Then $d_T(v_0)\in \{1,4\}$, $d_T(v_{t+1})\in \{1,4\}$, and $d_T(v_i)\in \{2,4\}$ for $i=2,\dots,t-1$.
Let $\{v_0,v_{t+1}\}$ have $r$ vertices with degree 1 and $s$ vertices with degree 4,
where $r+s=2$.
Let there are $x$ $(2,2)$-edges and $y$ $(2,4)$-edges of $T$ in $E(P)$,
where $x+y=t-1$.
Then $M_{22}\ge x$, $M_{24}\ge y+s$, and $M_{12}\ge r$.
By Lemma \ref{lem-0}, $M_{22}<M_{24}<M_{12}<0$. So
\begin{align*}
e^{SO_{red}}(T)\le &\frac{2n+2}{3}f(1,4)+\frac{n-5}{3}f(4,4)+rM_{12}+xM_{22}+(y+s)M_{24}.
\end{align*}

{\bf Subcase 2.1.}\ $t\ge 3$.

By Lemma \ref{lem-0}, $2M_{12}+2M_{24}<2M_{13}+M_{34}$. Then
\begin{align*}
e^{SO_{red}}(T)\le &\frac{2n+2}{3}f(1,4)+\frac{n-5}{3}f(4,4)+rM_{12}+xM_{24}+yM_{24}+sM_{12}
\\
=&\frac{2n+2}{3}f(1,4)+\frac{n-5}{3}f(4,4)+2M_{12}+(t-1)M_{24}
\\
\le &\frac{2n+2}{3}f(1,4)+\frac{n-5}{3}f(4,4)+2M_{12}+2M_{24}
\\
< &\frac{2n+2}{3}f(1,4)+\frac{n-5}{3}f(4,4)+2M_{13}+M_{34}
\\
=&2f(1,3)+f(3,4)+\frac{2n-5}{3}f(1,4)+\frac{n-7}{3}f(4,4).
\end{align*}

{\bf Subcase 2.2.}\ $t=2$.

Then $x=1$ and $y=0$. Since $n\ge 5$, we have $s\ge 1$ and $r\le 1$.
By Lemma \ref{lem-0}, $M_{12}+M_{22}+ M_{24}<2M_{13}+M_{34}$. So
\begin{align*}
e^{SO_{red}}(T)\le &\frac{2n+2}{3}f(1,4)+\frac{n-5}{3}f(4,4)+r M_{12}+M_{22}+s M_{24}
\\
\le &\frac{2n+2}{3}f(1,4)+\frac{n-5}{3}f(4,4)+M_{12}+M_{22}+ M_{24}
\\
< &\frac{2n+2}{3}f(1,4)+\frac{n-5}{3}f(4,4)+2M_{13}+M_{34}
\\
=&2f(1,3)+f(3,4)+\frac{2n-5}{3}f(1,4)+\frac{n-7}{3}f(4,4).
\end{align*}

{\bf Case 3.}\  $n_2(T)\ge 1$ and $n_3(T)\ge 1$.

{\bf Subcase 3.1.}\ $m_{23}(T)\ge 1$.

Note that $m_{12}(T)+m_{24}(T)\ge 1$, and $m_{13}(T)+m_{34}(T)\ge 2$.
By Lemma \ref{lem-0}, $M_{24}<M_{12}<0$, $M_{34}<M_{13}<0$, and $M_{23}<M_{34}<0$.
Then by (\ref{eq2.9}), we have
\begin{align*}
e^{SO_{red}}(T)=
&\frac{2n+2}{3}f(1,4)+\frac{n-5}{3}f(4,4)+M_{12}m_{12}(T)+M_{13}m_{13}(T)+M_{22}m_{22}(T)\\
&+M_{23}m_{23}(T)+M_{24}m_{24}(T)+M_{33}m_{33}(T)+M_{34}m_{34}(T)
\\
\le & \frac{2n+2}{3}f(1,4)+\frac{n-5}{3}f(4,4)\\
&+M_{23}+M_{12}(m_{12}(T)+m_{24}(T))+M_{13}(m_{13}(T)+m_{34}(T))
\\
\le & \frac{2n+2}{3}f(1,4)+\frac{n-5}{3}f(4,4)+M_{23}+M_{12}+2M_{13}
\\
<&\frac{2n+2}{3}f(1,4)+\frac{n-5}{3}f(4,4)+M_{34}+2M_{13}
\\
=&2f(1,3)+f(3,4)+\frac{2n-5}{3}f(1,4)+\frac{n-7}{3}f(4,4).
\end{align*}

{\bf Subcase 3.2}\ $m_{23}(T)=0$.

Then $m_{24}(T)\ge 1$, $m_{12}(T)+m_{24}(T)\ge 2$, $m_{34}(T)\ge 1$, and $m_{13}(T)+m_{34}(T)\ge 3$.
By Lemma \ref{lem-0}, $M_{22}<0$, $M_{33}<0$, $M_{24}<M_{12}<0$, and $M_{34}<M_{13}<0$.
Then by (\ref{eq2.9}), we have
\begin{align*}
e^{SO_{red}}(T)=
&\frac{2n+2}{3}f(1,4)+\frac{n-5}{3}f(4,4)+M_{12}m_{12}(T)+M_{13}m_{13}(T)+M_{22}m_{22}(T)\\
&+M_{24}m_{24}(T)+M_{33}m_{33}(T)+M_{34}m_{34}(T)
\\
\le &\frac{2n+2}{3}f(1,4)+\frac{n-5}{3}f(4,4)+M_{12}m_{12}(T)+M_{24}m_{24}(T)\\
&+M_{13}m_{13}(T)+M_{34}m_{34}(T)
\\
\le &\frac{2n+2}{3}f(1,4)+\frac{n-5}{3}f(4,4)+M_{12}+M_{24}+2M_{13}+M_{34}
\\
<&\frac{2n+2}{3}f(1,4)+\frac{n-5}{3}f(4,4)+M_{34}+2M_{13}
\\
=&2f(1,3)+f(3,4)+\frac{2n-5}{3}f(1,4)+\frac{n-7}{3}f(4,4).
\end{align*}

The lemma holds. \h

\begin{cor}\label{cor}\
Let $n\ge 5$ and $T\in {\cal CT}_n$.
If $n_2(T)+n_3(T)\ge 2$, then
\begin{align*}
e^{SO_{red}}(T)&<2f(1,3)+f(3,4)+\frac{2n-5}{3}f(1,4)+\frac{n-7}{3}f(4,4)
\\
&<f(1,2)+f(2,4)+\frac{2n-3}{3}f(1,4)+\frac{n-6}{3}f(4,4)
\\
&<\frac{2n+2}{3}f(1,4)+\frac{n-5}{3}f(4,4).
\end{align*}
\end{cor}

\p By Lemma \ref{lem-0}, we have $2M_{13}+M_{34}<M_{12}+M_{24}<0$. Then by Lemma \ref{lem-5},
\begin{align*}
e^{SO_{red}}(T)&<2f(1,3)+f(3,4)+\frac{2n-5}{3}f(1,4)+\frac{n-7}{3}f(4,4)
\\
&=\frac{2n+2}{3}f(1,4)+\frac{n-5}{3}f(4,4)+2M_{13}+M_{34}
\\
&<\frac{2n+2}{3}f(1,4)+\frac{n-5}{3}f(4,4)+M_{12}+M_{24}
\\
&=f(1,2)+f(2,4)+\frac{2n-3}{3}f(1,4)+\frac{n-6}{3}f(4,4)
\\
&<\frac{2n+2}{3}f(1,4)+\frac{n-5}{3}f(4,4).
\end{align*}
The result follows. \h

\section{Main results}
\hskip\parindent
Note that
\begin{align*}
\frac{2n+2}{3}f(1,4)+\frac{n-5}{3}f(4,4)=&\frac{2}{3}(n+1) e^3+\frac{1}{3}(n-5) e^{3 \sqrt{2}},
\\[2mm]
2f(1,3)+f(3,4)+\frac{2n-5}{3}f(1,4)+\frac{n-7}{3}f(4,4)=&2 e^2+e^{\sqrt{13}}+\frac{1}{3} (2 n-5) e^3
+\frac{1}{3} (n-7)e^{3 \sqrt{2}},
\end{align*}
\begin{align*}
f(1,2)+f(2,4)+\frac{2n-3}{3}f(1,4)+\frac{n-6}{3}f(4,4)=&e+e^{\sqrt{10}}+\frac{1}{3}(2n-3) e^3+\frac{1}{3} (n-6) e^{3 \sqrt{2}}.
\end{align*}
From Lemmas \ref{lem-1}--\ref{lem-5} and Corollary \ref{cor}, we have the following two theorems.

\begin{theorem}\label{thm-1}\
Let $n\ge 5$ and $T\in {\cal CT}_n$. Then
\begin{align*}
e^{SO_{red}}(T)\le \frac{2}{3}(n+1) e^3+\frac{1}{3}(n-5) e^{3 \sqrt{2}},
\end{align*}
the equality holds if and only if $n\equiv 2 (\!\!\!\!\mod 3)$, and $n_2(T)=n_3(T)=0$.
\end{theorem}

\begin{theorem}\label{thm-2}\
Let $n\ge 5$ and $T\in {\cal CT}_n$.

(1)\ If $n\equiv 2 (\!\!\!\!\mod 3)$, then
\begin{align*}
e^{SO_{red}}(T)\le \frac{2}{3}(n+1) e^3+\frac{1}{3}(n-5) e^{3 \sqrt{2}},
\end{align*}
the equality holds if and only if $n_2(T)=n_3(T)=0$.

(2)\ If $n\equiv 1 (\!\!\!\!\mod 3)$, then
\begin{align*}
e^{SO_{red}}(T)\le 2 e^2+e^{\sqrt{13}}+\frac{1}{3} (2 n-5) e^3+\frac{1}{3} (n-7)e^{3 \sqrt{2}},
\end{align*}
the equality holds if and only if $n_2(T)=0$, $n_3(T)=1$, $m_{13}(T)=2$, and $m_{34}(T)=1$.

(3)\ If  $n\equiv 0 (\!\!\!\!\mod 3)$, then
\begin{align*}
e^{SO_{red}}(T)\le e+e^{\sqrt{10}}+\frac{1}{3}(2n-3) e^3+\frac{1}{3} (n-6) e^{3 \sqrt{2}},
\end{align*}
the equality holds if and only if $n_2(T)=1$, $n_3(T)=0$, and $m_{12}(T)=m_{24}(T)=1$.
\end{theorem}

The results of Theorem \ref{thm-2} show that Conjecture \ref{conj-1} is negative.

\end{document}